\documentclass[12pt]{article}
\usepackage[cmex10]{amsmath}
\usepackage{amsfonts}
\usepackage{amsthm}
\usepackage{amssymb}
\usepackage{latexsym}
\usepackage[latin1]{inputenc}

\usepackage{amsmath}

\begin{document}

\newtheorem{theo}{Theorem}[section]
\newtheorem{lemma}[theo]{Lemma}
\newtheorem{cor}[theo]{Corollary}

\newtheorem{definition}[theo]{Definition}
\newtheorem{example}[theo]{Example}
\newtheorem{remark}[theo]{Remark}
\newtheorem{conj}[theo]{Conjecture}
\newtheorem{prop}[theo]{Proposition}
\newtheorem{reference}{}
\newcommand{\wh}{\widehat}
\newcommand{\ol}{\overline}

\newcommand{\Hom}{{\rm Hom}\,}           %  Homomorphism
\newcommand{\Aut}{{\rm Aut}}             %  Automorphism
\newcommand{\zzz}{\mathbb{Z}}        %  integer numbers
\newcommand{\car}{{\rm char}}
\def\F{{\mathbb F}}

\newcommand{\ndv}{ \ {\mid \kern -0.70 em {\scriptstyle \not}} \ \ }
%\def\lcm{\mathop{{\rm lcm}}}

\title{Group algebras and Coding Theory}

\author{Marinês Guerreiro       
\thanks{M. Guerreiro is with Departamento de Matemática, Universidade Federal de Vi\c cosa, CEP 36570-000 - Vi\c cosa-MG (Brasil), supported by 
CAPES, PROCAD 915/2010 and FAPEMIG, APQ CEX 00438-08 and PCE-00151-12(Brazil). E-mail: marines@ufv.br.}
\thanks{This work was presented at Groups, Rings and Group Rings 2014 on the occasion of the 70th Birthday of César Polcino Milies, Ubatuba-SP (Brasil) 21-26 July,  2014. }}

\maketitle

\begin{abstract}
%\boldmath
%Let $G$ be a finite group and $\F$ a field such that $\car(\F) \ndv |G|$. 
%Denote by 
%$\F G$ the group algebra of $G$ over $\F$.
%A (semisimple) abelian code is an ideal of $\F G$.  
%Two codes ${\cal I}_1$ and ${\cal I}_2$  of $\F G$ are 
%{\em $G$-equivalent} if there exists an automorphism $\psi$ of $G$ whose linear extension 
%to $\F G$ maps ${\cal I}_1$ onto ${\cal I}_2$.

Group algebras have been used in the context of Coding Theory since the beginning of the latter, but not in its full power. The work of Ferraz and Polcino Milies entitled {\it Idempotents in group algebras and minimal abelian codes (Finite Fields and their Applications, 13, (2007) 382-393)} gave origin to many thesis and papers linking these two subjects. In these works, the techniques of group algebras are mainly brought into play for the computing of the idempotents that generate the
 minimal codes and the minimum weight of such codes. 
In this paper I summarize the main results of the work done by doctorate students and research partners of Polcino Milies and Ferraz. 

\end{abstract}


\noindent {\bf Dedicated to Prof. C\'esar Polcino Milies on the occasion of his seventieth birthday.}

\maketitle

\section{Introduction}

The origins of Information Theory and Error Correcting Codes Theory are in the papers by 
 Shannon~\cite{shannon} and Hamming~\cite{hamming}, where they settled the theoretical foundations for 
 such theories. 

For a non empty finite set $A$, called {\bf alphabet}, a {\bf code} $C$ of {\bf length} $n$ is simply 
a proper subset of $A^n$ and an $n$-tuple $(a_0,a_1,\ldots, a_{n-1})\in C$ is called a {\bf word} of the code $C$. 

If $A=\F_q$ is a finite field with $q$ elements, then a {\bf linear code} $C$ of length $n$ is a proper subspace of $\F_q^n$. If $\dim C=k$ ($k<n$), then the number of words in $C$ is 
$q^k$. 

We shall call ``cyclic shift" the linear map $\pi : \F_q^n\longrightarrow\F_q^n$ such that $\pi (a_0, a_1, \ldots, a_{n-1}) = (a_{n-1}, a_0, a_1, \ldots, a_{n-2})$. 

A linear {\bf cyclic code} is a linear code $C$ that is invariant under the cyclic shift. This structure gives rise to fast-decoding algorithms, which is a considerable aspect regarding the conditions on communication.  

Consider the quotient ring $\displaystyle R_n = \frac{\F_q[x]}{<x^n-1>}$ and denote by $[f(x)]$ the class of the polynomial $f(x)$ in $R_n$. 
There is a natural  vector space isomorphism $\varphi : \F_q^n \longrightarrow R_n$ given by 
$$\varphi (a_0, a_1, \ldots, a_{n-1})= [a_0 + a_1 x+ \cdots + a_{n-1}x^{n-1}].$$  

Linear cyclic codes are often realized as ideals in $R_n$ and the cyclic shift is equivalent, via the isomorphism $\varphi$,
 to the multiplication 
by the class of $x$ in $R_n$. 

Group algebras may be defined in a more general setting, that is, for any group and over any field. However, we restrict the definitions and results below to finite groups and finite fields because this is the context for coding theory.

Let $G$ be a finite group written multiplicatively and $\F_q$ a finite field. 
The {\bf group algebra of $G$ over $\F_q$} is the set of all formal linear combinations
$$\alpha =\displaystyle \sum_{g \in G} \alpha_gg , \quad \mbox{where} \quad \alpha_g \in \F_q.$$

Given $\alpha = \displaystyle\sum_{g \in G} \alpha_gg$ and $\beta = \displaystyle\sum_{g \in G} \beta_gg$ we have
 $$\alpha = \beta \Longleftrightarrow \alpha_g = \beta_g, \ \ \mbox{ for all } g\in G.$$
 
 The support of an element $\alpha \in \F_q G$ is the set of elements of $G$ effectively appearing in $\alpha$; i.e.,
 $$ \operatorname{supp}(\alpha) = \{ g\in G \;|\; a_g\neq 0 \}.$$

 We define

 $$\left( \sum_{g \in G}\alpha_gg \right) + \left( \sum_{g \in G} \beta_gg \right) = \sum_{g \in G}(\alpha_g + \beta_g ) g.$$

$$\left( \sum_{g \in G}\alpha_gg \right)\left( \sum_{h \in G} \beta_h h \right) = \sum_{g,h \in G}(\alpha_g\beta_h)gh.$$

For $\lambda$ in $\F_q$, we define
$$\lambda \left( \sum_{g \in G}\alpha_gg \right) = \sum_{g \in G}(\lambda \alpha_g)g.$$

It is easy to see that, with the operations above, $\F_q G$ is an algebra over the field $\F_q$.

The {\bf weight} of an element $\alpha =\sum_{g\in G} a_g g\in \F_q G$ is the number of elements in its support; i.e.
$$w(\alpha) = |\{ g\;|\; a_g\neq 0 \}|.$$

For an ideal $I$ of $\F_q G$, we define the {\bf minimum weight} of $I$ as:
$$w(I) = min\{ w(\alpha)\;|\; \alpha\in I, \; \alpha\neq 0\}.$$

Let $C_n=\langle a\rangle$ denote a cyclic finite group of order $n$ generated by an element $a$. 
MacWilliams~\cite{mw1}  was the first one  to consider cyclic codes as ideals of the group ring 
$\F_q C_n$ which is easily proved to be isomorphic to $\frac{\F_q[x]}{<x^n-1>}.$  In $\F_q C_n$, the 
cyclic shift is equivalent to the multiplication of the elements of the code by $a$.

The following diagram helps us to understand the cyclic shift in these three different ways of considering a cyclic code.

$$ \begin{array} {cccrccrc}
 & C\subset \F_q^n & \stackrel{\varphi}{\longrightarrow}  & & R_n= \frac{\F_q[x]}{<x^n-1>} & \stackrel{\cong}{\longrightarrow} & & \F_q C_n=\F_q <a> \\
\mbox{cyclic} &    & & & & & & \\
 & \downarrow  & & & \bar{x} \;  \downarrow & & & a \; \downarrow \\
\mbox{shift} & & & & & &  & \\
 & C\subset \F_q^n & \stackrel{\varphi}{\longrightarrow}  & & R_n=\frac{\F_q[x]}{<x^n-1>} & \stackrel{\cong}{\longrightarrow} & & \F_q C_n =\F_q <a>
\end{array}$$

Extending these ideas, Berman~\cite{berman, berman1} and, independently, MacWilliams \cite{mw2} defined {\bf abelian codes} 
as ideals in finite abelian group algebras and, more generally, a {\bf group (left) code} was defined as an (left) ideal in a finite group algebra. Group codes were then   studied using 
ring and character-theoretical results. 

From now on, for a finite group $G$ and a finite field $\F_q$, we treat ideals in a group algebra $\F_q G$  as codes.
In this approach,  the length of the code is the order of the group $G$ and the dimension of a code $I$ is  its  dimension as an $\F_q$-subspace in $\F_q G$.  Length, dimension and minimum weight are the 
three parameters that define a linear code.

A group code is called {\bf minimal} if the corresponding ideal is minimal in the set of ideals of the group algebra. 
Keralev and Sol\'e in \cite{kersol} showed that many important codes can be realized as ideals in a group algebra, 
for example, the generalized Reed-Muller codes and generalized quadratic residue codes. These results are included in Section 9.1 of \cite{ker}. There is also a good treatment on the subject in~\cite{drenla}.

%{\bf S.D. Berman}, {\it Semisimple cyclic and abelian codes, II}. Kybernetika {\bf 3} (1967) 21-30.
 
%{\bf F.J.~MacWilliams}, {\it Binary codes which are ideals in the group algebra of an abelian group}. Bell System Tech. %Journal {\bf 44} (1970) 987-1011.

A word of warning is necessary here, because the expression ``group code" may  also have some other meanings. For example,  in Computer Science, sometimes  group codes consist of $n$ linear block codes which are subgroups of $G^n$, where $G$ is a finite abelian group, as in~\cite{BE,FT}.

Usually in the papers that present techniques to compute the idempotents that generate the codes, 
character theory is used in the context of polynomials, as it can be seen in~\cite{arorap1, arorap2, bakshiraka1, bakshiraka, barasha,  miller, poli1, poli2, bdrs}. Sometimes the expressions for the idempotents are not very ``reader friendly". 
Moreover, the character theory and polynomial approaches in the computation of idempotents 
 did not fully  explore the structure of the group underneath the group algebra that defines the 
 underlying set for the codes. 
 
\

Summarizing the work of 
 Ferraz and Polcino Milies~\cite{FM}, in Section~\ref{section2}
we define the idempotents using  subgroups of the group and establish the basic theorems that are used in this  
work. We emphasize that in~\cite{FM}, they gave simpler proofs for computing idempotents, dimension and minimum weight for minimal cyclic and abelian codes of length $2^kp^n$, generalizing and comparing their results with the ones in~\cite{arorap1, arorap2}. Full details of Ferraz and Polcino Milies´ paper are described in the Master´s Dissertation of 
Luchetta~\cite{valeria}.

In Section~\ref{sec3} we discuss some topics of non-abelian group codes, including some equivalence questions. In Section~\ref{sec4} we describe some results on codes of length $2^n$, for a natural $n\geq 1$ and in Section~\ref{sectionequiv} we treat some aspects of equivalence of abelian codes. In Section~\ref{sec6} we summarize some results on cyclic and abelian codes of length $p^nq^m$, for $p$ and $q$ distinct primes and $n, m \geq 1$ and in Section~\ref{sec7} we explore some facts on codes over rings. 

\section{Subgroups and idempotents}\label{section2}

We recall that an element in the group algebra $\F_q G$ is called {\bf central} if it commutes with every other element of the algebra. A non-zero central idempotent $e$ is called {\bf primitive} if it cannot be decomposed in the form $e=e'+e''$, where $e'$ and $e''$ are both non-zero central idempotents such that $e'e''=e''e'=0$. For $\operatorname{char}(\F_q)  \ndv |G|$, the  group algebra $\F_q G$ is semisimple and the primitive central idempotents are the generators of the minimal two-sided ideals. Two idempotents $e',e''$ are  {\bf  orthogonal} if $e'e''=e''e'=0$.

The primitive central idempotents of the rational group algebra $\mathbb{Q}G$ were 
computed in~\cite[Theorem~VII.1.4]{goodcesar} in the case $G$ abelian; in \cite[Theorem 2.1]{jlp} 
when $G$ is nilpotent; in~\cite[Theorem 4.4]{OdRS} in a more general context and 
in~\cite[Theorem 7]{BdR} an algorithm to compute the primitive idempotents is given. 

In what follows, we shall establish a correspondence between primitive idempotents of $\F_q G$ and certain subgroups of an abelian group $G$.

Let $G$ be a finite (abelian)  group and $\F_q$  a field such that  $\car(\F_q)\ndv |G|$.
Given a subgroup $H$ of $G$, denote
\begin{equation} \label{idempotentH}
\widehat{H} = \frac{1}{|H|}\sum_{h\in H}h
\end{equation}  which is an idempotent of 
$\F_q G$ and, for an element $x\in G$, set $\wh x = \wh{\left< x\right>}$.

It is known that the idempotent $\wh{G}$ is always primitive (see \cite[Proposition 3.6.7]{MS}).

\begin{definition} \label{subgrpcocyl}
Let $G$ be an abelian group. A subgroup $H$  of $G$ is called
 a \textbf{co-cyclic subgroup} if the factor group $G/H\neq \{ 1\}$ is  cyclic.
\end{definition}
We use the notation
$$
{\mathcal{S}}_\mathrm{cc}(G)=\{H\, |\, H \mbox{ is a co-cyclic subgroup of } G\}.
$$

For a finite group $G$, denote by $\operatorname{exp}(G)$ the {\bf exponent} of $G$ which is the smallest positive integer $t$ such that $g^t=1$, for all $g\in G$. A group $G$ is called a {\bf $p$-group} if its exponent is a power of a given prime $p$. In particular, this means that the order of every element of $G$ is itself a power of $p$.

Let $G$ be a finite abelian  $p$-group and $\F_q$  a field such that  $\car(\F_q)\ndv |G|$.
For each co-cyclic subgroup $H$ of $G$, we can construct an idempotent of $\F_q G$.
In fact, we remark that, since $G/H$ is a cyclic $p$-group, there exists a unique subgroup $H^{\sharp}$ of $G$ containing $H$ such that $|H^{\sharp}/H|=p$.
Then $e_H=\wh H - \wh{H^{\sharp}}$ is an idempotent and we consider the set 
\begin{equation}\label{idempabelian}
\{\wh G\}\cup \{ e_H=\wh H - \wh{H^{\sharp}}\,|\,H\in {\mathcal{S}}_\mathrm{cc}(G)\}.
\end{equation}

We recall the following results that are used throughout this paper.

In the case of a rational abelian group algebra $\mathbb{Q} G$, the set~\eqref{idempabelian} is the set  of all primitive central idempotents~\cite[Theorem~1.4]{goodcesar}. \\

\begin{theo}{\em \cite[Lemma~5]{FM}} \label{3.44}
Let $p$ be a prime  
integer and $G$ a finite abelian group of exponent $p^n$ and $\F_q$ a finite field with $q$ elements such that $p \ndv q$. Then~\eqref{idempabelian} 
  is a set of pairwise orthogonal idempotents of $\F_q G$ whose sum is equal to $1$, i.e., 
  \begin{equation}
1= \wh G\,+\,\sum_{H\in {\mathcal{S}}_\mathrm{cc}(G)} e_H,
\end{equation}
where $1$ also denotes the identity  element in $\F_q G$.
\end{theo}

In our next statement, we denote by $U(\zzz_{p^n})$ the set of invertible elements of the ring $\zzz_{p^n}$ of integers modulo $p^n$; $\bar{q}$ denotes
the class of the integer $q$ in $\zzz_{p^n}$ and, when it is invertible, $o(\bar{q})$ denotes its multiplicative order; i.e., the least positive integer $m$ such that $\bar{q}^m=\bar{1}$.

\begin{theo}{\em \cite[Theorem~4.1]{FM}}  \label{3.4} Under the same hypotheses of Theorem~\ref{3.44}, 
the set~\eqref{idempabelian} is the set of all primitive idempotents of $\F_q G$ if and only  if    $o(\bar{q}) = \phi (p^n)$ in $U(\zzz_{p^n})$,
  where $\phi$ denotes Euler's totient function.
\end{theo}

For positive integers $r$ and $m$, we shall denote by $\bar{r}\in \mathbb{Z}_m$ the image of $r$ in the ring of integers modulo $m$. Then, for an element $g$ in a group $G$, define $ \mathcal{G}_{g}=\{ g^r\,|\,\gcd (r, o(g))=1\}= \{ g^r\,|\, \bar{r}\in U(\mathbb{Z}_{o(g)}) \}$.
The following theorem gives us conditions on the exponent $e$ of the group $G$ and the size $q$ of the finite field that 
satisfy Theorem~\ref{3.4}.

\begin{cor}{\em\cite[Teorema 7.10]{PauloMartin}}\label{corolarioprincipal}
Let $\mathbb{F}_{q} $ be a finite field with $q$ elements and $G$ a finite abelian group with exponent $e$ such that $\operatorname{gcd}\left(q,|G|\right)=1$. Then $C_{g}=\mathcal{G}_{g}$, for all $g \in G$, if only if one of following conditions holds, where $\phi$ denotes Euler´s totient function:
\begin{enumerate}
\item[(a)] $e=2$ and $q$ is odd;
\item[(b)]$e=4$ and $q \equiv 3 (\operatorname{mod}\,\,4)$;
\item[(c)]$e=p^{n}$ and $o(q)= \phi(p^{n})$ in $U\left(\mathbb{Z}_{p^{n}}\right)$;
\item[(d)]$e=2p^{n}$ and $o(q)= \phi(p^{n})$ in $U\left(\mathbb{Z}_{2p^{n}}\right)$.
\end{enumerate}
\end{cor}

\begin{theo} {\em\cite[Lemma 3]{FM}} \label{osidempotentesprimitivosparaumpgrupo}
Let $G=\langle g\rangle$ be a cyclic group with order $p^{n}$ and $\mathbb{F}_{q} $ a finite field with $q$ elements such that $\overline{q}$ generates $U\left(\mathbb{Z}_{p^{n}}\right)$. Consider
$$
G=G_{0} \supset G_{1} \supset ... \supset G_{n}=\{1\}
$$
the descending chain of all subgroups of $G$. Then a complete set of primitive idempotents in $\mathbb{F}_{q} G$ is:
\begin{equation}\label{idempotPrimitivos}
e_{0}=\widehat{G}=\frac{1}{p^{n}} \sum_{g \in G}g \quad and \quad e_{i}=\widehat{G_i}-\widehat{G_{i-1}}, \mbox{ for } 1 \leq i \leq n,
\end{equation}
with $G_i=<g^{p^i}>$, for $1 \leq i \leq n$.
\end{theo}

As the authors comment in~\cite{FM}, a straightforward computation shows that these
are the same idempotents given in~\cite[Theorem 3.5]{arorap2}, though there they are expressed in
terms of cyclotomic cosets.

The idempotent generators of minimal ideals in the case of cyclic groups of order $2p^n$
now follow easily from the previous results.
\begin{theo} {\em\cite[Theorem 2.6]{arorap2}} 
Let $\F_q$ be a finite field with $q$ elements and
$G$ a cyclic group of order $2p^n$, $p$ an odd prime, such that $o(\bar{q}) = \phi(p^n)$ in $U(\mathbb{Z}_{2p^n} )$.
Write $G = C \times  A$, where $A$ is the $p$-Sylow subgroup of $G$ and $C = \{ 1, t\}$ is its $2$-Sylow
subgroup. If $e_i$, for $ 0\leq  i\leq n$, denote the primitive idempotents of $\F_q A$, then the primitive
idempotents of $\F_q G$ are
\begin{equation}
\frac{1+t}{2}\,e_i,\quad \frac{1-t}{2}\,e_i,\quad 0\leq i\leq n.
\end{equation}
\end{theo}

More generally,
\begin{theo}{\em \cite[Theorem 4.2]{FM}}
Let $p$ be an odd prime,  $A$ be an abelian $p$-group of exponent $2p^r$ and $\F_q$ be a finite field with $q$ 
elements such that $o(q)= \phi(p^{r})$ in $U\left(\mathbb{Z}_{2p^{r}}\right)$.
Write $A = E × B$, where $E$ is an elementary abelian $2$-group and $B$ a $p$-group. Then the
primitive idempotents of $\F_q A$ are products of the form $e\cdot f$ , where $e$ is a primitive idempotent
of $\F_q E$ and $f$ a primitive idempotent of $\F_q B$.
\end{theo}

Section~5 of~\cite{FM} is devoted to the computation of dimension and minimum weight of codes generated by
the idempotents presented in previous theorems. 
For non-cyclic abelian groups, we may also apply the ideas above to construct idempotents. 
In~\cite{FGMequiv}, the following results are presented in details.

For a finite abelian group $G$, we write $G=G_{p_1}\times \cdots \times G_{p_t}$, where
$G_{p_i}$ denotes the $p_i$-Sylow subgroup of $G$, for the distinct 
prime numbers $p_1,\ldots, p_t$.

\begin{lemma}{\em \cite[Lemma II.5]{FGMequiv}}\label{co-cyclic}
Let $G=G_{p_1}\times  \cdots \times G_{p_t}$ be a finite abelian group and $H\in {\mathcal{S}}_\mathrm{cc}(G)$.
Write $H=H_{p_1}\times  \cdots \times H_{p_t}$, where $H_{p_i}$ is the $p_i$-Sylow subgroup of $H$. Then 
each subgroup $H_{p_i}$  is co-cyclic in $G_{p_i}$, $1\leq i\leq t$.
\end{lemma}

With the notation above, for each $H\in {\mathcal{S}}_\mathrm{cc}(G)$, define an idempotent $e_H\in \F_q G$ as follows. For each $1\leq i\leq t$,
either $H_{p_i}=G_{p_i}$ or there exists a unique subgroup $H_{p_i}^{\sharp}$ such that $[H_{p_i}^{\sharp}:H_{p_i}]=p_i$. Thus, let
$e_{H_{p_i}}=\wh{G_{p_i}}$ or $e_{H_{p_i}}=\wh{H_{p_i}}-\wh{H_{p_i}^{\sharp}}$, respectively, and define
\begin{equation}\label{idempH}
e_H = e_{H_{p_1}}e_{H_{p_2}}\cdots e_{H_{p_t}}.
\end{equation}

For any other $K\in {\mathcal{S}}_\mathrm{cc}(G)$, with $K\neq H$, we have $K_{p_i}\neq H_{p_i}$, for some $1\leq i\leq t$, and, by Theorem~\ref{3.44},
$e_{H_{p_i}}e_{K_{p_i}}=0$, hence $e_He_K=0$. It is easy to see that  $\wh{G} e_H =0$, for all $H\in {\mathcal{S}}_\mathrm{cc}(G)$.

 Thus, we have the following.

\begin{prop}{\em \cite[Proposition II.6]{FGMequiv}}\label{baseB}
Let $G$ be a finite abelian group and $\F_q$ a field such that $\car(\F_q)\ndv |G|$. 
Then 
\begin{equation}\label{basebeta}
{\mathcal{B}}=\{ e_H\,|\,H \in {\mathcal{S}}_\mathrm{cc}(G) \} \cup \{ \widehat{G} \} 
\end{equation}
is a set of orthogonal idempotents of $\F_q G$, where $e_H$ is defined as in~\eqref{idempH}.
\end{prop}

A similar construction of idempotents for rational group algebras of abelian groups is given in~\cite[Section VII.1]{goodcesar}. 
For the rational case, these idempotents are primitive while for finite fields this is usually not true.

Now, we extend Theorem~\ref{3.44} to finite abelian groups.

\begin{lemma}{\em \cite[Lemma II.7]{FGMequiv}}\label{sum=1}
Let $G$ be a finite abelian group and $\F_q$ a field such that $\car(\F_q)\ndv |G|$. 
Then, in the group algebra $\F_q G$, we have
\begin{equation}\label{soma=1}
1= \wh G\,+\,\sum_{H\in {\mathcal{S}}_\mathrm{cc}(G)} e_H.
\end{equation}
\end{lemma}
The following lemma starts the discussion about the relation between idempotents and certain subgroups of the abelian group, which we elaborate in more details in Section~\ref{sectionequiv}.
\begin{lemma}{\em \cite[Lemma II.8]{FGMequiv}}\label{uniqueco-cyc}
Let $G$ be a finite abelian group and $\F_q$ a field such that $\car(\F_q) \ndv |G|$. For each primitive idempotent $e\in\F_q G$, $e\neq \wh{G}$,
there
exists a unique $H\in {\mathcal{S}}_\mathrm{cc}(G)$ such that $e\cdot e_H = e$. Also, $e\cdot e_K=0$, for any other $K\in {\mathcal{S}}_\mathrm{cc}(G)$.%, with $K\neq H$.
\end{lemma}

\section{Non Abelian Codes}\label{sec3}

\subsection{Dihedral and Quaternion Codes} 

As a natural way to proceed, Theorem~\ref{3.44} is used by Dutra~\cite{dutra, FDM} in her Ph.D. thesis 
to compute idempotents for non abelian group codes, particularly, for dihedral and quaternion groups. 

For $n\geq 1$, Dutra considered the semisimple group algebras $\F_qD_n$ of the
dihedral groups $D_n = \langle a, b\,|\,a^n=b^2=1,\,bab=a^{-1} \rangle$ 
over a finite field $\F_q$ and gave conditions under which the
number of its simple components is minimum, that is, the same as for the rational group algebra
$\mathbb{Q} D_n$. These conditions are stated in the following theorem.

\begin{theo}{\em\cite[Teorema 2.2]{dutra}}\label{teo1Dutra}
Let $\F_q$ be a field with $q$ elements and $D_n$ the dihedral group with $2n$
elements such that $\gcd (q, 2n) = 1$. Let $p, p_1$ and $p_2$ be 
distinct odd primes and $m, m_1$ and $m_2$ be positive integers. Then
$ \F_q D_n$ and  $\mathbb{Q} D_n$ have the same number of simple components if and only if one of the following conditions occurs: 

(i) $n = 2$ or $4$ and $q$ is odd.

(ii) $n = 2^m$, with  $m \geq 3$ and  congruent to $3$ or $5$ modulo $8$.

(iii) $n = p^m$ and the class $\bar{q}$ generates the group of units $U(\mathbb{Z}_{p^m})$.

(iv) $n = p^m$, the class $\bar{q}$ generates the group $U^2(\mathbb{Z}_{p^m})=\{ x^2\,|\, x\in U(\mathbb{Z}_{p^m})\}$ and 
$-1$ is not a square modulo $p^m$.

(v) $n = 2p^m$ and the class $\bar{q}$ generates the group of units $U(\mathbb{Z}_{p^m})$.

(vi) $n = 2p^m$,  $\bar{q}$ generates the group $U^2(\mathbb{Z}_{p^m})=\{ x^2\,|\, x\in U(\mathbb{Z}_{p^m})\}$ and 
$-1$ is not a square modulo $p^m$.

(vii) $n = 4p^m$, $4$ divides $\phi(p^m)$ and the class $\bar{q}$ generates the group $U(\mathbb{Z}_{p^m})$.

(viii) $n = 4p^m$, $4$ does not divide $\phi(p^m)$, $q\equiv 1 (\operatorname{mod }4)$ and the class $\bar{q}$ generates the group $U(\mathbb{Z}_{p^m})$.

(ix) $n = 4p^m$, $4$ does not divide $\phi(p^m)$, $q\equiv -1 (\operatorname{mod } 4)$ and the class $\bar{q}$ has order $\phi(p^m)/2$.

(x) $n = p_1^{m_1} p_2^{m_2}$, with $\gcd (\phi(p_1^{m_1}), \phi(p_2^{m_2} )) = 2$ and $q$ or $-q$ has order $\phi(p_1^{m_1} p_2^{m_2})/2$ modulo $p_1^{m_1} p_2^{m_2}$.

(xi) $n = 2 p_1^{m_1} p_2^{m_2}$, with $\gcd (\phi(p_1^{m_1}), \phi(p_2^{m_2} )) = 2$ and $q$ or $-q$ has order $\phi(p_1^{m_1} p_2^{m_2})/2$ modulo $p_1^{m_1} p_2^{m_2}$.
\end{theo}

Under such conditions, Dutra
computed the set of minimal codes of $\F_qD_n$,  
their dimensions, minimum weights and bases for these codes as follows.

\begin{theo} {\em \cite[Proposição~3.1]{dutra}}
Let $q$ and $n$ be integers related as in conditions (i) and (ii) of Theorem~\ref{teo1Dutra}.
If $\mathcal{C} $ is a dihedral code of length $2n$ generated by the idempotent $e$, then $\mathcal{C}$ has dimension and
minimum weight described in the table below.
$$
\begin{array}{|c|c|c|}
\hline
e & \operatorname{dim}_{\F_q\mathcal{C}} & w(\mathcal{C} )\\  \hline 
\widehat{b}\widehat{a} & 1 & 2^{m+1}\\  \hline
(1 - \widehat{b})\widehat{a} & 1 &  2^{m+1} \\ \hline 
\widehat{b}(\widehat{a^2} - \widehat{a}) & 1 &  2^{m+1} \\ \hline
(1-\widehat{b})(\widehat{a^2} - \widehat{a}) & 1 &  2^{m+1} \\ \hline
(\widehat{a^{2^{i}}} - \widehat{a^{2^{i-1}}}) & 2^i &  2^{m-i+1} \\ \hline
\end{array}
$$

\end{theo}

\begin{theo} {\em\cite[Proposição~3.2]{dutra}}
Let $q$ and $n$ be integers related as in conditions (iii) and (iv) of Theorem~\ref{teo1Dutra}.
If $\mathcal{C} $ is a dihedral code of length $2n$ generated by the idempotent $e$, then $\mathcal{C}$ has dimension and
minimum weight described in the table below.
$$
\begin{array}{|c|c|c|}
\hline
e & \operatorname{dim}_{\F_q\mathcal{C}} & w(\mathcal{C} )  \\ \hline
\widehat{b}\widehat{a} & 1 & 2p^{m} \\ \hline
(1 - \widehat{b})\widehat{a} & 1 &  2p^{m} \\ \hline 
(\widehat{a^{p^{i}}} - \widehat{a^{p^{i-1}}}) & 2\phi(p^i) &  2p^{m-i} \\ \hline
\end{array}
$$

\end{theo}

\begin{theo} {\em \cite[Proposição~3.3]{dutra}}
Let $q$ and $n$ be integers related as in conditions (v) to (ix) of Theorem~\ref{teo1Dutra}.
For $n=p_1^{m_1}p_2^{m_2}$ with $p_1=2, m_1=1$ or $2$ and $p_2$ an odd prime, if $\mathcal{C} $ 
is a dihedral code of length $2n$ generated by the idempotent $e_1e_2$, then $\mathcal{C}$ has dimension and
minimum weight described in the table below.
$$
\begin{array}{|c|c|c|c|}
\hline
e_1 & e_2 & \operatorname{dim}_{\F_q\mathcal{C}} & w(\mathcal{C} )\\  \hline 
\widehat{b}\widehat{C_{p_1^{m_1}}} & \widehat{C_{p_2^{m_2}}} & 1 & 2p_1^{m_1}p_2^{m_2}\\  \hline
(1-\widehat{b})\widehat{C_{p_1^{m_1}}} & \widehat{C_{p_2^{m_2}}} & 1 & 2p_1^{m_1}p_2^{m_2}\\  \hline
\widehat{b}(\widehat{C_{p_1^{m_1-1}}}-\widehat{C_{p_1^{m_1}}}) & \widehat{C_{p_2^{m_2}}} & 1 & 2p_1^{m_1}p_2^{m_2}\\  \hline
(1-\widehat{b})(\widehat{C_{p_1^{m_1-1}}}-\widehat{C_{p_1^{m_1}}}) & \widehat{C_{p_2^{m_2}}} & 1 & 2p_1^{m_1}p_2^{m_2}\\  \hline
\widehat{C_{p_1^{m_1}}} & \widehat{C_{p_2^{m_2-j}}}-\widehat{C_{p_2^{m_2-j+1}}} & 2\phi (p_2^j) & 2p_1^{m_1}p_2^{m_2-j}\\
  \hline
\widehat{C_{p_1^{m_1-i}}}-\widehat{C_{p_1^{m_1-i+1}}} & \widehat{C_{p_2^{m_2}}} & 2\phi (p_1^j) & 2p_1^{m_1-i}p_2^{m_2}\\  \hline
\widehat{C_{p_1^{m_1-i}}}-\widehat{C_{p_1^{m_1-i+1}}} & \widehat{C_{p_2^{m_2-j}}}-\widehat{C_{p_2^{m_2-j+1}}}
 & 2\phi (p_1^j)2\phi (p_2^j) & 4p_1^{m_1-i}p_2^{m_2-j}\\  \hline
\end{array}
$$
\end{theo}

\begin{theo} {\em\cite[Proposição~3.4]{dutra}}
Let $q$ and $n=p_1^{m_1}p_2^{m_2}$, with $p_1$ and $p_2$ odd distinct prime numbers, integers related as in condition (x) 
of Theorem~\ref{teo1Dutra}. If $\mathcal{C} $ 
is a dihedral code of length $2n$ generated by the idempotent $e_1e_2$, then $\mathcal{C}$ has dimension and
minimum weight described in the table below.
$$
\begin{array}{|c|c|c|c|}
\hline
e_1 & e_2 & \operatorname{dim}_{\F_q\mathcal{C}} & w(\mathcal{C} )\\  \hline 
\widehat{b}\widehat{C_{p_1^{m_1}}} & \widehat{C_{p_2^{m_2}}} & 1 & 2p_1^{m_1}p_2^{m_2}\\  \hline
(1-\widehat{b})\widehat{C_{p_1^{m_1}}} & \widehat{C_{p_2^{m_2}}} & 1 & 2p_1^{m_1}p_2^{m_2}\\  \hline
\widehat{C_{p_1^{m_1}}} & \widehat{C_{p_2^{m_2-j}}}-\widehat{C_{p_2^{m_2-j+1}}} & 2\phi (p_2^j) & 2p_1^{m_1}p_2^{m_2-j}\\
  \hline
\widehat{C_{p_1^{m_1-i}}}-\widehat{C_{p_1^{m_1-i+1}}} & \widehat{C_{p_2^{m_2}}} & 2\phi (p_1^j) & 2p_1^{m_1-i}p_2^{m_2}\\  \hline
\widehat{C_{p_1^{m_1-i}}}-\widehat{C_{p_1^{m_1-i+1}}} & \widehat{C_{p_2^{m_2-j}}}-\widehat{C_{p_2^{m_2-j+1}}}
 & 2\phi (p_1^j)2\phi (p_2^j) & 4p_1^{m_1-i}p_2^{m_2-j}\\  \hline
\end{array}
$$
\end{theo}

\begin{theo} {\em\cite[Proposição~3.5]{dutra}}
Let $q$ and $n=2p_1^{m_1}p_2^{m_2}$, with $p_1$ and $p_2$ odd distinct prime numbers, integers related as in condition (xi) 
of Theorem~\ref{teo1Dutra}. If $\mathcal{C} $ 
is a dihedral code of length $2n$ generated by the idempotent $e_0e_1e_2$, then $\mathcal{C}$ has dimension and
minimum weight described in the table below.
{\tiny
$$
\begin{array}{|c|c|c|c|c|}
\hline
e_0 & e_1 & e_2 & \operatorname{dim}_{\F_q\mathcal{C}} & w(\mathcal{C} )\\  \hline 
\widehat{b}\widehat{C_2} & \widehat{C_{p_1^{m_1}}} & \widehat{C_{p_2^{m_2}}} & 1 & 4p_1^{m_1}p_2^{m_2}\\  \hline
(1-\widehat{b})\widehat{C_2} & \widehat{C_{p_1^{m_1}}} & \widehat{C_{p_2^{m_2}}} & 1 & 4p_1^{m_1}p_2^{m_2}\\  \hline
\widehat{b}(1-\widehat{C_2}) & \widehat{C_{p_1^{m_1}}} & \widehat{C_{p_2^{m_2}}} & 1 & 4p_1^{m_1}p_2^{m_2}\\  \hline
(1-\widehat{b})(1-\widehat{C_2}) & \widehat{C_{p_1^{m_1}}} & \widehat{C_{p_2^{m_2}}} & 1 & 4p_1^{m_1}p_2^{m_2}\\  \hline
\widehat{C_2} & \widehat{C_{p_1^{m_1}}} & \widehat{C_{p_2^{m_2-j}}}-\widehat{C_{p_2^{m_2-j+1}}} & 2\phi (p_2^j) & 4p_1^{m_1}p_2^{m_2-j}\\
  \hline
\widehat{C_2} & \widehat{C_{p_1^{m_1-i}}}-\widehat{C_{p_1^{m_1-i+1}}} & \widehat{C_{p_2^{m_2}}} & 2\phi (p_1^j) & 4p_1^{m_1-i}p_2^{m_2}\\  \hline
(1-\widehat{C_2}) & \widehat{C_{p_1^{m_1}}} & \widehat{C_{p_2^{m_2-j}}}-\widehat{C_{p_2^{m_2-j+1}}} & 2\phi (p_2^j) & 4p_1^{m_1}p_2^{m_2-j}\\
  \hline
(1-\widehat{C_2}) & \widehat{C_{p_1^{m_1-i}}}-\widehat{C_{p_1^{m_1-i+1}}} & \widehat{C_{p_2^{m_2}}} & 2\phi (p_1^j) & 4p_1^{m_1-i}p_2^{m_2}\\  \hline
\widehat{C_2} & \widehat{C_{p_1^{m_1-i}}}-\widehat{C_{p_1^{m_1-i+1}}} & \widehat{C_{p_2^{m_2-j}}}-\widehat{C_{p_2^{m_2-j+1}}}
 & 2\phi (p_1^j)2\phi (p_2^j) & 8p_1^{m_1-i}p_2^{m_2-j}\\  \hline
(1-\widehat{C_2}) & \widehat{C_{p_1^{m_1-i}}}-\widehat{C_{p_1^{m_1-i+1}}} & \widehat{C_{p_2^{m_2-j}}}-\widehat{C_{p_2^{m_2-j+1}}}
 & 2\phi (p_1^j)2\phi (p_2^j) & 8p_1^{m_1-i}p_2^{m_2-j}\\  \hline
\end{array}
$$
} %close small letter
\end{theo}

Similar results were obtained by Dutra~\cite[Capítulos 4 e 5]{dutra} for group codes over the quaternion groups.

\subsection{Metacyclic Codes and Equivalence Questions}\label{metaequiv}

A group $G$ is {\bf metacyclic} if  it contains a normal cyclic subgroup $H$ such that $G/H$ is also cyclic. It is easy to prove that a finite metacyclic group has the following presentation
\begin{equation}
G=\left\langle a,\,b\,|\, a^m=1,\,b^n=a^s,\, bab^{-1}=a^i\right\rangle,
\end{equation}
with $a$ and $b$ such that $H=\langle a\rangle$ and $G/H=\langle bH\rangle$, for $m, n\in \mathbb{N}$ and $1\leq s, i\leq m$ such that $s|m,\;m|s(i-1),\, i<m, \,\gcd(i,m)=1$. 
For $s=m$, we say that $G$ is a {\bf split metacyclic group} and, in this case, $G$ is the semi-direct product $G=\langle a\rangle\rtimes\langle b\rangle$.

Earlier approaches on non-abelian metacyclic codes include results obtained by Sabin~\cite{sabin} and Sabin and 
Lomonaco~\cite{sabinLom}, where we also find the following definition of equivalence of codes.

\begin{definition}\label{combequiv}
Let $G$ and $H$ be two finite groups of the same order and $\F_q$ a field. A {\bf combinatorial equivalence} 
is a vector space isomorphism $\psi: \F_q G \longrightarrow \F_q H$ induced by a bijection $\psi : G\longrightarrow H$.
\vspace{0.2cm}
 
Two codes $C\subset \F_q G$ and $\tilde{C}\subset \F_q H$ are said to be {\bf combinatorially equivalent} if there exists a
combinatorial equivalence $\psi: \F_q G \longrightarrow \F_q H$ such that $\psi(C)=\tilde{C}$. 
\end{definition}

For $G$ a metacyclic finite group such that $\gcd (q, |G|)=1$, Sabin and Lomonaco~\cite{sabinLom}, by
using group representation theory, 
proved that codes generated by central idempotents in $\F_q G$ are combinatorially equivalent to abelian codes. This motivated the search for left minimal codes in $\F_q G$.

Considering the group algebra of a non-abelian split metacyclic group $G$ over a finite field $\F_q$, Assuena~\cite{samirtese} in his Ph.D. thesis, found a necessary condition under which $\F_q G$ has the minimum number of simple components. 
\begin{theo}{\em\cite[Teorema 2.1.16]{samirtese}}
Let $G$ be a metacyclic group and $\F_q$ a finite field with $q$ elements such that $\gcd (q, |G|)=1$. If the number of simple components of the group algebra $\F_q G$ is minimal, then $U(\mathbb{Z}_n)=\langle \bar{q}\rangle$ and 
$U(\mathbb{Z}_m)=\langle \bar{i}\rangle \langle \bar{q}\rangle$.
\end{theo}

In his thesis, Assuena used the structure of the group to determine the minimal metacyclic codes for a
non abelian split metacyclic group of order $p^m\ell^n$, with $p$ and $\ell$ odd prime numbers, under the conditions
that $\F_q G$ is semisimple and the number of simple components of $\F_q G$ is minimum.

For $D_{p^m}$, the dihedral group of order $2p^m$, and $\F_q$ a finite field such that $\gcd (q,2p^m)=1$, 
he constructs left minimal codes that {\bf are not} combinatorially equivalent to abelian codes and also exhibits
one case where a left minimal code is more efficient then the abelian ones of the same length, giving a positive 
answer to a conjecture of Sabin and Lomonaco~\cite{sabinLom}.

 Further studies on group codes are given in~\cite{BdRS}, where it is defined a {\bf (left) $G$-code} as any linear 
 (left) code of length $n$ over a field $\F_q$ which is the image of a (left) ideal of a group algebra via an 
 isomorphism $\F_q G\longrightarrow \F_q^n$ which maps the finite group 
 $G$ of order $n$ to the standard basis of $\F_q^n$.  Their ideas are used in~\cite{schafer} to study two-sided and abelian group ring codes and in~\cite{PGMMN1}, where García Pillado et al. first communicated an example of a non-abelian $S_4$-code over $\F_5$. The full proof of this computacional construction was given later in~\cite{PGMMN2}. New examples of non-abelian $G$-codes are given in~\cite{PGMMN3} and, particularly, using the group $SL(2;\F_3)$ instead of the
symmetric group, they prove, without using a computer for it, that there is a code
over $\F_2$ of length $24$, dimension $6$ and minimal weight $10$. This code has greater minimum
distance than any abelian group code having the same length and dimension over $\F_2$,
and, moreover, it has the greatest minimum weight among all binary linear
codes with the same length and dimension.
 
 In~\cite{eliaG} Elia and García Pillado give an overview of the properties of
 ideal group codes defined as principal ideals
in the group algebra of a finite group $G$ over a finite field $\F_q$  and present their encoding and
syndrome decoding. They also describe in detail a correction of a single
error, using syndromes. 
 
\section{Cyclic codes of length $2^m$}\label{sec4}

Codes are usually considered over the binary field $\F_2$. For cyclic codes of length $2^m$, with a natural $m\geq 1$, 
over a field of odd size, the results obtained using a polynomial approach  by Bakshi and Raka~\cite{bakshiraka1}, Pruthi~\cite{pruthi}, Sharma et al.~\cite{sbdr}, Sharma et al.~\cite{sbr} and using the group algebra approach, by Prado~\cite{prado} in her Ph.D. thesis, are
essentially the same. In Chapter 2, Prado states the general facts:

\begin{theo}{\em \cite[Lema 2.1.1]{prado}}\label{janete1}
Let $G=\langle a\rangle$ be a finite cyclic group of order $2^m, m\geq 1$ and $\F_q$ a finite field of odd characteristic. Let
$$
G=G_0 \supset G_1 \supset \cdots \supset G_m=\{ 1\}
$$
be the descending chain of all subgroups of $G$, with $G_i = \langle a^{2^i}\rangle$ and $|G_i|=2^{m-i}$. Then the elements
$e_0 = \widehat{G}$ and $e_i=\widehat{G_i}-\widehat{G_{i-1}}$, with $1\leq i\leq m$, form a set of orthogonal idempotents of $\F_q G$ such that $e_0+e_1+\cdots +e_m=1$.
\end{theo}

\begin{theo}{\em\cite[Lema 2.2.1]{prado}} Under the same hypothesis of Theorem~\ref{janete1},
let $I_i=\F_q G e_i$, with  $1\leq i\leq m$, be the ideals of $\F_q G$ generated by the idempotents $e_i$ of 
 Theorem~\ref{janete1}. Then 
\begin{eqnarray}
\operatorname{dim}(I_0) & = & 1,\qquad d(I_0)=|G|=2^m \nonumber\\
\operatorname{dim}(I_i) & = & 2^{i-1},\quad d(I_i)=|G|=2^{m-i+1}, \mbox{ for } 1\leq i\leq m.  \nonumber
\end{eqnarray} 
\end{theo}

The notion of a visible code was given by Ward~\cite{ward}, where he defines a {\bf visible basis} for a code as 
a basis where all its elements have the same weight. Prado also proved the following for codes of length $2^m$.

\begin{theo} {\em \cite[Proposição 2.3.1]{prado}} Under the same hypothesis of Theorem~\ref{janete1},  for $1\leq i\leq m$, the 
set 
$$
B_i = \{ e_i, ae_i, a^2e_i, \ldots, a^{2^{i-1}-1}e_i\}
$$
is a visible basis for the code $I_i=\F_q G e_i$.
\end{theo}

In her thesis~\cite{prado}, Prado  studied in details the minimal codes generated by primitive idempotents in $\F_q C_{2^m}$, 
with $q$ odd. She
considered four cases: $q\equiv 1 (\operatorname{mod} 8)$, $q\equiv 3 (\operatorname{mod} 8)$, 
$q\equiv 5 (\operatorname{mod} 8)$ and $q\equiv 7 (\operatorname{mod} 8)$. The order of 
$q (\operatorname{mod} 2^m)$, the number of simple components of $\F_q C_{2^m}$ and the computation of idempotents 
are different for each one of these cases. For $q\equiv 3 (\operatorname{mod} 8)$ and 
$q\equiv 5 (\operatorname{mod} 8)$ a complete discussion is presented in the thesis and the other cases are exemplified with particular examples. Here is the case $q\equiv 3 (\operatorname{mod} 8).$ 

\begin{theo}{\em\cite[Proposição 3.1.1]{prado}}
Let $\F_q$ be a field with $q$ elements such that $q\equiv 3 (\operatorname{mod} 8)$ and $G=\langle a\,|\, a^{2^m}=1\rangle$ 
be a cyclic group of order $2^m$. The following elements of the group algebra $\F_q G$
\begin{center}

$ e_0=\displaystyle\frac{1+a+a^2+\cdots + a^{2^m-1}}{2^m}$

$ e_1=\displaystyle\frac{1-a+a^2-\cdots - a^{2^m-1}}{2^m}$

$ e_2=\displaystyle\frac{1-a^2+a^4-\cdots - a^{2^m-2}}{2^{m-1}}$

$ e_3=(1-a^4)\displaystyle\frac{(1+a^{2^3}+\cdots + a^{2^m-2^3})(2+\alpha a + \alpha a^3)}{2^m}$

$ e'_3=(1-a^4)\displaystyle\frac{(1+a^{2^3}+\cdots + a^{2^m-2^3})(2-\alpha a - \alpha a^3)}{2^m}$

$ e_4=(1-a^8)\displaystyle\frac{(1+a^{2^4}+\cdots + a^{2^m-2^4})(2+\alpha a^2 + \alpha a^{3\cdot 2})}{2^{m-1}}$

$ e'_4=(1-a^8)\displaystyle\frac{(1+a^{2^4}+\cdots + a^{2^m-2^4})(2-\alpha a^2 - \alpha a^{3\cdot 2})}{2^{m-1}}$

\

$\ldots ,$

\

$ e_{m-1}=(1-a^{2^{m-2}})\displaystyle\frac{(1+a^{2^m-1})(2+\alpha a^{2^{m-4}} + \alpha a^{3\cdot 2^{m-4}})}{2^{4}}$

$ e'_{m-1}=(1-a^{2^{m-2}})\displaystyle\frac{(1+a^{2^m-1})(2-\alpha a^{2^{m-4}} - \alpha a^{3\cdot 2^{m-4}})}{2^{4}}$

$ e_{m}=(1-a^{2^{m-1}})\displaystyle\frac{(2+\alpha a^{2^{m-3}} + \alpha a^{3\cdot 2^{m-3}})}{2^{3}}$

$ e'_{m}=(1-a^{2^{m-1}})\displaystyle\frac{(2-\alpha a^{2^{m-3}} - \alpha a^{3\cdot 2^{m-3}})}{2^{3}}$
\end{center}
form a complete set of primitive idempotents of $\F_q G$, with $\alpha^2=-2$ in $\F_q$.
\end{theo}

In Chapter 4 of her thesis, Prado simplifies results of Poli~\cite{poli1} in order to obtain a clearer description 
of the principal nilpotent ideals of a group algebra of finite abelian groups in a modular case (i.e., when $\operatorname{char} (\F_q)$ divides the order of the group $G$). She also exemplifies the process of lifting idempotents modulo a nilpotent ideal.

\section{More on equivalence of abelian codes}\label{sectionequiv}

The question of equivalence in Coding Theory has many approaches. In Section~\ref{metaequiv}, we have defined combinatorial equivalence. In~\cite{miller}, we found the following for abelian codes. Here $G$ stands for a finite abelian group and
$\F_q $ is a finite field with $q$ elements.

\begin{definition} Two abelian codes ${\cal I}_1$ and ${\cal I}_2$  are {\em $G$-equivalent} if there exists an automorphism $\theta$ of $G$ whose linear extension to $\F_qG$ maps ${\cal I}_1$ on ${\cal I}_2$.
\end{definition}

The following statements also appeared in~\cite{miller}.

\

\noindent{\bf Theorem A}~\cite[Theorem 3.6]{miller}  {\it Let  $G$ be a finite abelian group of odd 
order and exponent $n$ and denote by
$\tau(n)$ the number of divisors of $n$. Then there exist  precisely $\tau(n)$ non $G$-equivalent minimal abelian codes 
 in $\F_2 G$}.

\

\noindent{\bf Theorem B}~\cite[Theorem 3.9]{miller}  {\it Let  $G$ be a finite abelian group of odd order. Then two minimal abelian codes in $\F_2 G$
 are $G$-equivalent if and only if they have the same weight distribution.}

\

Unfortunately both statements are not correct. The errors arise from the assumption, implicit in the last paragraph of~\cite[p. 167]{miller}, that  if  $e$ and $f$ are primitive idempotents of $\F_2 C_m$ and $\F_2 C_n$, respectively, then $ef$ is a primitive idempotent of
$\F_2[C_m\times C_n]$. To the best of our knowledge, these results have not been used in a wrong way in the literature.

We first communicated the following counterexamples to both Theorems A and B in~\cite{FGM}.

\begin{prop}{\em\cite[Proposition 3.1]{FGM}} \label{cp2cp}
Let $p$ be an odd prime such that $\bar 2$ generates $U(\zzz_{p^2})$ and $G=\left< a\right>\times \left< b\right>$ an abelian group, with  $o(a)=p^2$ and $o(b)=p$. Then $\F_2G$ has four inequivalent minimal codes, namely, the ones generated by the idempotents
$e_0=\wh G$, $e_1= \widehat{b} - \widehat{\left< a^p\right>\times \left< b\right>} $,
$e_2= \widehat{a} - \wh G$ and $e_3= \widehat{\left< a^p\right>\times \left< b\right>} - \wh G.$

Also all minimal codes  of $\F_2 G$ are described in Table~\ref{tabela1} with their dimension and weight.
\begin{table}[h]

{\small
$$\begin{array}{|l|c|c|c|} \hline \label{cp2cptable}
\mbox{\rm Code} & \mbox{\small\rm Primitive Idempotent} & \mbox{\rm Dimension} & \mbox{\small\rm Minimum Weight}  \\
             \hline\hline
I_0 & e_0=\wh{a}\wh{b} = \wh G       &  1      & p^3  \\ \hline
I_{1} & e_1= \wh{b} - \wh{\left< a^p\right>\times \left< b\right>}       &  p^2-p      & 2p   \\ \hline
I_{1j} & e_{1j}= \wh{a^{jp}b} - \wh{\left< a^p\right>\times \left< b\right>}       &  p^2-p      & 2p  \\
&   j=1,\ldots, p-1 & & \\ \hline
I_{2} & e_2=\wh{a} - \wh{G}       &  p-1      & 2p^2   \\ \hline
I_{2i} & e_{2i}=\wh{ab^i} - \wh{G}       &  p-1      & 2p^2  \\
& i=1,\ldots, p-1 & & \\ \hline
I_{3} & e_3=\wh{\left< a^p\right>\times \left< b\right>} - \wh G        & p-1      & 2p^2  \\ \hline
    \hline
\end{array}
 $$
} % end of small command
\caption{Minimal codes in $\F_2(C_{p^2}\times C_p)$\label{tabela1}} 
\end{table}
Moreover, the minimal inequivalent codes $I_2$ and $I_3$ have the same weight distribution.
\end{prop}

In~\cite[Proposition 4.2]{FGM} we showed that Theorem~A holds in the special case of 
minimal codes in $\F_2(C_{p^n}\times C_{p^n})$ and, in~\cite[Theorem V.3]{FGMequiv}, we generalize these
result for $G$ a direct product of $m\geq 2$ copies of a cyclic group $ C_{p^n}$, as follows.

\begin{prop}{\em \cite[Proposition V.3]{FGMequiv}}\label{propcprgeral}
Let $m$ and $r$ be positive integers and $p$ a prime number. If $G=(C_{p^r})^{m}$  is a finite abelian $p$-group
and $\F_q$ is a field of $ \car (\F_q)\neq p$.  Then a primitive idempotent of $\F_q G$, different from $\wh{G}$, is of the form
$\wh K\cdot e_h$, where $K$ is a subgroup of $G$ isomorphic to
$(C_{p^r})^{m-1}$ and $e_h$ is a primitive
idempotent of $\F_q\langle h\rangle$, where $h\in G$ is such that $G=\langle h\rangle \times K$ and $\langle h\rangle \cong C_{p^r}$.
\end{prop}

This result can be applied as follows.

\begin{cor} {\em \cite[Corollary V.4]{FGMequiv}}
Let $m$ and $r$ be positive integers, $p$ a prime number, a finite abelian $p$-group 
$G= (C_{p^r})^{m}$  and $\F_q$
a finite field with $q$ elements such that  $o(\bar{q})=\phi(p^r)$ in $U(\zzz_{p^r})$. Then
 the minimal abelian codes  in $\F_q G$ are as follows, where $h$ and $K$ are as   in Proposition~\ref{propcprgeral}.
$$
\begin{array}{|c|c|c|}
\hline\hline
\mbox{\rm Primitive Idempotent} &  \mbox{\rm Dimension} & \mbox{\rm Weight}\\ \hline\hline
\wh G     & 1  & p^{rm}\\ \hline
 \wh{K}(\wh{h^{p}} - \wh{h}) & p-1 &  2p^{r(m-1)+(r-1)}\\ \hline
 \wh{K}(\wh{h^{p^2}} - \wh{h^{p}}) & p(p-1) &  2p^{r(m-1)+(r-2)}\\ \hline
 \wh{K} (\wh{h^{p^3}} - \wh{h^{p^2}}) & p^2(p-1) & 2p^{r(m-1)-(r-3)} \\ \hline
\cdots  & \cdots &  \\ \hline
\wh{K}(\wh{h^{p^{i}}} - \wh{h^{p^{i-1}}}) & p^{i-1}(p-1) & 2p^{r(m-1)-(r-i)}\\ \hline
\cdots  & \cdots & \\ \hline
 \wh{K} (1 - \wh{h^{p^{r-1}}}) & p^{r-1}(p-1) & 2p^{r(m-1)} \\
 \hline\hline
\end{array}
$$

Consequently, the number of non $G$-equivalent minimal abelian codes is $r+1=\tau(p^r)$.
\end{cor}

\begin{cor}{\em \cite[Corollary V.5]{FGMequiv}}
Let $n,m\geq 2$ be   integers, $G=(C_n)^{m}$ an abelian group and $\F_q$ a finite field
such that $\gcd ( q, n) = 1$.
Then the primitive idempotents of $\F_qG$ are of the form $\wh K\cdot e_h$, where $K$ is a subgroup of $G$ isomorphic to
$ (C_n)^{m-1}$, $h\in G$ is such that $G=K\times \langle h\rangle$ and $e_h$ is a primitive
idempotent of $\F_q\langle h\rangle$.%, for $1\leq i\leq \tau(n)$.
\end{cor}

\begin{theo} {\em \cite[Theorem V.6]{FGMequiv}}
Let $G = C^n$ be a direct product
of cyclic groups isomorphic to one another, of exponent $n$, and $\F_q$ a finite field such that $char(\F_q)\ndv  |G|$.
Then, the number of non $G$-equivalent minimal abelian codes is precisely $\tau(n)$. 
\end{theo}

\

We fully discussed the $G$-equivalence of abelian codes and established in~\cite[Section III]{FGMequiv} a relation
between the classes of equivalence of $G$-equivalent codes and some classes of isomorphisms of subgroups of $G$, as follows.

We say that two subgroups $H$ and $K$ of a group $G$ are {\bf\em $G$-isomorphic} if there exists an automorphism $\psi\in \Aut(G)$ such that $\psi(H)=K$.

Notice that isomorphic subgroups are not necessarily $G$-isomorphic. For example, for a prime $p$,
if $G=\langle a\rangle\times \langle b\rangle$ with $o(a)=p^2$ and $o(b)=p$, then $\langle a^p\rangle$ and $\langle b \rangle$ are isomorphic, as they are both cyclic groups of order $p$. However, they are not $G$-isomorphic, since $\langle b\rangle$ is contained, as a subgroup of index $p$, only in $\langle a^p\rangle\times \langle b\rangle$ while
$\langle a^p\rangle$ is contained in $\langle a\rangle$ and in $\langle a^ib\rangle$, for all $1\leq i\leq p-1$. An automorphism of $G$ carrying one to the other would preserve also inclusions.

We shall denote by
$ \mathcal{P}(\F_q G)$  the set of all primitive idempotents of $\F_q G$. Recall the notion of co-cyclic subgroup (Definition~\ref{subgrpcocyl}). Then, under the same hypotheses of Lemma~\ref{uniqueco-cyc}, the following map is well-defined
\begin{equation}\label{mapPhi}
\begin{array}{clc}
  \Phi\; : \;\mathcal{P}(\F_q G) & \longrightarrow & {\mathcal{S}}_\mathrm{cc}(G) \cup \{G\}\\
 e\neq \wh{G} &\longmapsto & \Phi(e)=H_e,\\
 \wh{G} & \longmapsto & G
\end{array}
\end{equation}
where $H_e$ is the unique co-cyclic subgroup of $G$ such that $e\cdot e_{H_e}=e$. 

\begin{theo}{\em \cite[Theorem II.9]{FGMequiv}}\label{somaeH}
Let $G$ be a finite abelian group, $\F_q$ a field such that $\car(\F_q) \ndv |G|$ and $H\in {\mathcal{S}}_\mathrm{cc}(G)$.
Then $e_H$ is the sum of all primitive idempotents $e\in \mathcal{P}(\F_q G) $ such that $\Phi(e)=H$.
\end{theo}

The study of the $G$-equivalence of ideals involves to know how the group of automorphisms $\Aut(G)$ acts on the
lattice of the subgroups of $G$ and hence on the idempotents in the group algebra which arise from these subgroups.
From now on, we use the same notation for an automorphism of the group $G$ and its linear extension to
the group algebra $\F_q G$. The following results from~\cite{FGMequiv} relate subgroups in $G$ and idempotents in $\F_qG$.

\begin{lemma}{\em\cite[Lemma III.1]{FGMequiv}}\label{invaut}
Let $G$ be a finite abelian group, $H\in {\mathcal{S}}_\mathrm{cc}(G)$ and $e_H$ its corresponding idempotent defined as in~\eqref{idempH}. Then, for any $\psi\in \Aut(G)$,  we have $\psi(e_H)=e_{\psi(H)}$ and
$\psi(\widehat{G})=\widehat{G}$. 
 \end{lemma}

For finite abelian groups, Propositions~\ref{prop1},~\ref{prop2} and~\ref{recprop1} below establish a correspondence between
$G$-equivalent minimal ideals in $\F_q G$  and $G$-isomorphic subgroups of $G$.

\begin{prop}{\em\cite[Proposition III.2]{FGMequiv}}\label{prop1}
Let $G$ be a finite abelian group and $\F_q$ a field such that $\car(\F_q)\ndv |G|$.
If $e,\,e'\in \mathcal{P}(\F_q G)$ are such that $\psi(e)=e'$, for some automorphism $\psi\in \Aut(G)$ linearly extended to $\F_q G$, then
$$ \psi(H_e)=H_{\psi(e)}=H_{e'},$$
i.e., $H_e$ and $H_{e'}$ are $G$-isomorphic.
\end{prop}

 We set
$
{\mathcal{L}}\Aut(G) = \{ \psi\in\Aut(G)\,|\,\psi(H)=H, \mbox{ for all } H\leq G \}.
$

\begin{prop}{\em\cite[Proposition III.7]{FGMequiv}}\label{prop2}
Let $G$ be a finite abelian group and $\F_q$ a field such that $\car(\F_q)\ndv |G|$.
If $e',\,e''\in \mathcal{P}(\F_q G)$ are both different from $\wh{G}$ and $H_{e'}=H_{e''}$, then there exists an automorphism
$\psi\in \mathcal{L}\Aut(G)$ whose linear extension to $\F_q G$ maps $e'$ to $e''$.
\end{prop}

The following is the converse of Proposition~\ref{prop1}.

\begin{prop}{\em\cite[Proposition III.8]{FGMequiv}}\label{recprop1}
Let $G$ be a finite abelian group and $\F_q$ a field such that $\car(\F_q)\ndv |G|$.
If \ $e',\,e''\in \mathcal{P}(\F_q G)$, both different from $\wh{G}$, are such that $\psi(H_{e'})=H_{e''}$, for some $\psi\in \Aut(G)$, then
there exists an automorphism $\theta\in \Aut(G)$ whose linear extension to $\F_q G$ maps  $e'$ to $e''$, i.e., the ideals 
of $\F_q G$  generated by
$e'$ and $e''$ are $G$-equivalent.
\end{prop}

As an application of Propositions~\ref{prop1} and~\ref{recprop1}, in~\cite[Section IV]{FGMequiv}
 we consider the minimal codes in $\F_2(C_{p^n}\times C_{p})$, for an odd prime $p$ and $n\geq 3$.
Its proof is similar to the proof of Proposition~\ref{cp2cp}.
This gives a whole family of counterexamples to Theorem~A.

\begin{prop}{\em\cite[Theorem IV.3]{FGMequiv}}\label{codescpncp}
Let $n\geq 3$ be a positive integer and $p$ an odd prime such that $\bar 2$ generates
$U(\zzz_{p^n})$ and $G=\left< a\right>\times \left< b\right>$ be an abelian group, with $o(a)=p^n$ and $o(b)=p$.
Then the minimal codes of $\F_2 G$ are described in Table~\ref{tabela2}. Moreover, 
there are $2n$ inequivalent minimal codes in $\F_2(C_{p^n}\times C_p)$.
\begin{table}[h]
$$
\begin{array}{|l|c|c|} \hline \label{cpncptable}
\mbox{\rm Code} & \mbox{\rm Dimension} & \mbox{\rm Weight}  \\
       \hline\hline
I_0 =\langle\wh{a}\wh{b}\rangle = \langle\wh G\rangle       &  1      & p^{n+1} \\ \hline
I_{1} =\langle\wh{\left< a^p\right>\times \left< b\right>} - \wh{G} \rangle     &  p-1      & 2p^n   \\ \hline
I_{1i} =\langle \wh{ab^i} - \wh{G}\rangle & p-1 & 2p^n \\
 i=0,\ldots, p-1 & &  \\ \hline
I_{2} =\langle\wh{\left< a^{p^2}\right>\times \left< b\right>} -\wh{\left< a^p\right>\times \left< b\right>}\rangle &  p(p-1) & 2p^{n-1}   \\ \hline
I_{2i} =\langle\wh{a^{p}b^i} - \wh{\left< a^p\right>\times \left< b\right>}    \rangle &  p(p-1)      & 2p^{n-1} \\
 i=1,\ldots, p-1 &  & \\ \hline
%I_{3} & e_3=\wh{\left< a^{p^3}\right>\times \left< b\right>} -   \wh{\left< a^{p^2}\right>\times \left< b\right>} & p^2(p-1) & 2p^{n-2} & \\ \hline
%I_{3i} & e_{3i}=\wh{a^{p^2}b^i} - \wh{\left< a^{p^2}\right>\times \left< b\right>}     &  p^2(p-1)      & 2p^{n-2}  &  i=1,\ldots, p-1\\ \hline
\ldots     & \ldots & \\ \hline
I_{k} =\langle\wh{\left< a^{p^k}\right>\times \left< b\right>} -   \wh{\left< a^{p^{k-1}}\right>\times \left< b\right>}\rangle & p^{k-1}(p-1) & 2p^{n-k+1}  \\ \hline
I_{ki} =\langle\wh{a^{p^{k-1}}b^i} - \wh{\left< a^{p^{k-1}}\right>\times \left< b\right>}\rangle     &  p^{k-1}(p-1)      & 2p^{n-k+1}  \\
i=1,\ldots, p-1 & & \\ \hline
\ldots &     \ldots  & \\ \hline
I_{n-1} =\langle\wh{\left< b\right>} -   \wh{\left< a^{p^{n-2}}\right>\times \left< b\right>}\rangle & p^{(n-1)}(p-1) & 2p  \\ \hline
I_{n-1,i} =\langle\wh{a^{p^{n-1}}b^i} - \wh{\left< a^{p^{n-2}}\right>\times \left< b\right>}  \rangle   &  p^{(n-1)}(p-1)      & 2p \\
  i=1,\ldots, p-1 & & \\ \hline
\hline
\end{array}
$$
\caption{Minimal codes in $\F_2(C_{p^n}\times C_p)$ \label{tabela2}}
\end{table}
\end{prop}

In the first column of Table~\ref{tabela3} we give a complete list of representatives of   classes of $G$-isomorphisms of subgroups of $C_{p^n}\times C_{p}$ and, in the second column, we list the corresponding representatives of
   $G$-equivalent classes of minimal codes of $\F_2(C_{p^n}\times C_{p})$. 
\begin{table}[h]
$$
\begin{array}{|c|c|} \hline \label{cpncpgrp}
\mbox{\rm Subgroups } & \mbox{\rm Codes}  \\ \hline\hline
G & I_0 = \langle \wh G \rangle \\ \hline
\left< a \right> & I_{11}= \langle \wh a  - \wh{G} \rangle \\ \hline
\left< a^p \right> \times \left< b \right> & I_{1}= \langle \wh{\left< a^p \right> \times \left< b \right>} - \wh{G} \rangle \\ \hline
\left< a^pb \right> &  I_{21} =\langle \wh{a^pb}-\wh{\left< a^p\right>\times \left< b\right>}  \rangle \\ \hline
\left< a^{p^2} \right> \times \left< b \right> &  I_{2} =\langle \wh{\left< a^{p^2} \right> \times \left< b \right>}-
\wh{\left< a^p\right>\times \left< b\right>}  \rangle \\ \hline
\ldots & \ldots \\ \hline
\left< a^{p^k}b \right> &  I_{k+1,1} =\langle \wh{a^{p^k}b}-\wh{\left< a^{p^k}\right>\times \left< b\right>}  \rangle \\ \hline
\left< a^{p^{k+1}} \right> \times \left< b \right> &  I_{k+1} =\langle \wh{\left< a^{p{k+1}} \right> \times \left< b \right>}-
\wh{\left< a^{p^k}\right>\times \left< b\right>}  \rangle \\ \hline
\ldots & \ldots \\ \hline
\left< b \right> &  I_{n-1} =\langle \wh b-
\wh{\left< a^{p^n{-1}}\right>\times \left< b\right>}  \rangle \\ \hline
   \hline
\end{array}
$$
\caption{ \label{tabela3}}
\end{table}

\section{Cyclic and abelian codes of length $p^nq^m$} \label{sec6}

\subsection{Binary abelian codes}

In~\cite{CFGM}, we considered finite abelian groups of type $G=G_p\times G_q$, for distinct primes $p$ and $q$ such that 
$G_p$ is a $p$-group, $G_q$ is a $q$-group satisfying 
the following conditions which will allow us to use the results in ~\cite{FM}:
\begin{eqnarray} \label{hypopq}
(i)& & \gcd (p-1,\,q-1 ) =  2, \nonumber \\
 (ii) & &  \  \bar{2}\;
\mbox{{\em generates the groups of units}}\; U(\zzz_{p^2})  \mbox{ {\em and} } 
\;U(\zzz_{q^2}) \\
(iii) & & \gcd(p-1,q) = \gcd(p,q-1)=1. \nonumber
\end{eqnarray}

The hypothesis $(i)$ above implies that at least one of the primes $p$ and $q$ is congruent to 3 (mod 4). 
In this section, to fix notations, we shall always assume that $q\equiv 3 \ $ (mod 4). As a code (ideal) generated by a primitive idempotent is isomorphic to a field, condition $(i)$ also helps us to have some control on the number 
of simple components that appear in the group algebra $\F_q (G_p\times G_q)$, because of the following elementary facts of Number Theory.

\begin{lemma}\label{lemars}
Let $\ell$ be a positive prime number and $\,r,\,s\in \mathbb{N}^*$. Then
$$\,\F_{\ell^r}\otimes_{\F_{\ell}} \F_{\ell^s}\,\cong
\,\gcd ( r,s)\cdot \F_{\ell^{\operatorname{lcm} (r,s)}}.$$
\end{lemma}

\begin{lemma}\label{lemae1e2}
Let $r,\,s\in{\mathbb N}$ be non-zero elements such that $\gcd(r,s)=2$. Let $u\in \F_{2^{r}}\,$ and $v\in
\F_{2^{s}}\,$ be elements  satisfying the equation $\,x^2+x+1=0$.
Then
\begin{equation}\label{compuv}
\F_{2^{r}}\otimes_{\F_{2}}\,\F_{2^{s}}\cong \F_{2^{\frac{rs}{2}}}\oplus\,\F_{2^{\frac{rs}{2}}}
\end{equation}
and $\,e_1\,=\,(u\otimes v)\,+\,(u^2\otimes v^2)\,$
and $\,e_2\,=\,(u\otimes v^2)\,+\,(u^2\otimes v)\,$ are the primitive idempotents generating to the simple components of~\eqref{compuv}.
\end{lemma}

Methods to determine idempotent generators for minimal cyclic codes were given in \cite{bakshiraka, barasha, singhpruthi} using representation theory. We develop our results without appealing to representation theory, working inside the group algebra. 

%For each subgroup $H$ in $G_p$ such that $G_p/H\neq 1$ is cyclic, consider the idempotent $e_H =\wh{H} -\wh{H^*}$ as above %and, similarly, consider primitive idempotents of the form $e_K = \wh{K} -\wh{K^*}$ for $G_q$.

For two co-cyclic subgroups $H$ of $G_p$ and $K$ of $G_q$, consider the respective idempotents $e_H =\wh{H} -\wh{H^*}$ in $\F_2 G_p$  and $e_K = \wh{K} -\wh{K^*}$ in $\F_2 G_q$.
Clearly $\wh{G_p}\cdot \wh{G_q}=\wh{G_p\times G_q}$ is a primitive idempotent of $\F_2 G=\F_2 (G_p\times G_q)$.  

It is ease to prove that  idempotents of the form $\wh{G_p} \cdot e_K$  and $e_H\cdot \wh{G_q}$ are primitive in $\F_2 G$. 
We proved that each idempotent of the form $e_H\cdot e_K$ decomposes as the sum of two primitive idempotents in $\F_2 G$,
by the following argument. For $e_H=\wh{H}-\wh{H^*}$, set $a\in H^*\setminus H$ (hence $aH$ is a generator of $H^*/H$). Set

\begin{equation}\label{uf2p}
u\,=\,\left\{
\begin{array}{ll}
a^{2^0}\,+\,a^{2^2}\,+\,\cdots\,+\,a^{2^{p-3}}, &
\mbox{if} \;p\equiv 1 (\operatorname{mod} 4)\;\mbox{or}\\
1\,+\,a^{2^0}\,+\,a^{2^2}\,+\,\cdots\,+\,a^{2^{p-3}},
& \mbox{if} \;p\equiv 3 (\operatorname{mod} 4) \end{array}\right.
\end{equation}
and
\begin{equation}\label{vf2p}
u'\,=\,\left\{
\begin{array}{ll}
a^{2}\,+\,a^{2^3}\,+\,\cdots\,+\,a^{2^{p-2}}, &
\mbox{if} \;p\equiv 1 (\operatorname{mod} 4)\;\mbox{or}\\
1\,+\,a^{2}\,+\,a^{2^3}\,+\,\cdots\,+\,a^{2^{p-2}},
& \mbox{if} \;p\equiv 3 (\operatorname{mod} 4) \end{array}\right.
\end{equation}

For $e_K=\wh{K}-\wh{K^*}$, set $b\in K^*\setminus K$ and define $v$ and $v'$ as in~\eqref{uf2p} and~\eqref{vf2p} replacing $a$ by $b$. As $\gcd(p^{r-1}(p-1),q^{s-1}(q-1))=2$ we can apply Lemma~\ref{lemae1e2} to see that 
\begin{eqnarray*}
e_1(H,K) & = & u\wh{H}\cdot v\wh{K} \;+\; u'\wh{H}\cdot v'\wh{K} \mbox{ and}\\
e_2(H,K) & = & u\wh{H}\cdot v'\wh{K} \;+\; u'\wh{H}\cdot v\wh{K}
\end{eqnarray*}
are primitive orthogonal idempotents such that $e_1+e_2 = e_He_K$.

Hence, we have shown the following.

\begin{theo}{\em\cite[Theorem III.1]{CFGM}}\label{teorema}
Let $G_p$ and $G_q$ be abelian $p$ and $q$-groups, respectively satisfying the conditions in \eqref{hypopq}. For a group $G$, denote by $S(G)$ the set of subgroups $N$ of $G$ such that $G/N\neq 1$ is cyclic. Then
 the set of primitive idempotents in $\F_2 [G_p\times G_q]$ is:
\begin{eqnarray*}
\wh{G_p} \cdot \wh{G_q}, && \\
  \wh{G_p}\cdot e_K,& &  K\in S(G_q),\\
 e_H\cdot \wh{G_q},& &  H\in S(G_p),\\
 e_1(H,K), \ e_2(H,K),& &   H\in S(G_p), K\in S(G_q).
 \end{eqnarray*}
\end{theo}

Particularly, in~\cite[Section IV]{CFGM} we compute, for  each minimal code of $\,\F_2(C_p\times C_q)$,  the  generating primitive idempotent, its dimension and give explicitly a basis for it over $\F_2$. In~\cite[Theorem IV.7]{CFGM} we presented the results on minimum weight for  these codes. In~\cite[Theorem V.1]{CFGM} we deal with the case $\F_2(C_{p^m}\times C_{q^n})$, for $m\geq 2, n\geq 2$ and extend this technique for three primes as follows.

\begin{theo}{\em\cite[Theorem IV.10]{CFGM}}\label{3primes}
Let $\,p_1,\,p_2\,$ and $\,p_3\,$ be three distinct positive odd prime numbers such that
$\,\gcd (p_i-1,\,p_j-1)=2$, for $1\leq i\neq j\leq 3$, and $\,\bar{2}$
generates the groups of units $U(\zzz_{p_i})$. Then the primitive idempotents of
the group algebra $\F_2G$ for the finite abelian group $\,G=C_{p_1}\times C_{p_2}\times C_{p_3}$, with $C_{p_1}=<a>$, $C_{p_2}=<b>$ and
$C_{p_2}=<c>$, are

$e_0=\hat{a}\hat{b}\hat{c}$, 
$e_1=\hat{a}\hat{b}(1-\hat{c})$, 
$e_2=\hat{a}(1-\hat{b})\hat{c}$,
$e_3=(1-\hat{a})\hat{b}\hat{c}$,

$e_4=(uv+u^2v^2)\hat{c}$,
$e_5=(u^2v+uv^2)\hat{c}$
$e_6=(uw+u^2w^2)\hat{b}$, 

$e_7=(u^2w+uw^2)\hat{b}$
$e_8=(vw+v^2w^2)\hat{a}$,
$e_9=(v^2w+vw^2)\hat{a}$

$e_{10}=(1-\hat{a})(1-\hat{b})(1-\hat{c})+ u^2v^2w+uvw^2$

$e_{11}=(1-\hat{a})(1-\hat{b})(1-\hat{c})+ u^2v^2w^2+uvw$

$e_{12}=(1-\hat{a})(1-\hat{b})(1-\hat{c})+ u^2vw+uv^2w^2$ and

$e_{13}=(1-\hat{a})(1-\hat{b})(1-\hat{c})+ uv^2w+u^2vw^2$,\\
where $u=u(a),\, v=v(b),\, w=w(c)$ are defined as in \eqref{uf2p}. 
\end{theo}

Comparing and using both the group algebra techniques of~\cite{CFGM, FM} with the polynomial techniques of~\cite{bakshiraka}, Bastos and Guerreiro~\cite{gusCNMAC, gusICMCTA} improved the presentation of minimal idempotents of length $p^nq$ given in~\cite{kumararora}, correcting some coefficients in their expressions.

\subsection{Codes of length $p^n$ also for non-cyclic abelian groups}

Let $\F_q$ be a finite field with $q$ elements and $G$ a cyclic group of order $p^n$ generated by $a$ such that $\gcd (q,p)=1$. 
Then the group algebra $\F_q G$ is semisimple and each of its ideals is a direct sum of minimal ones. Under the 
conditions (b) and (c) of Corollary~\ref{corolarioprincipal}, the minimal ideals (codes) 
are generated by the primitive idempotents given by Theorem~\ref{osidempotentesprimitivosparaumpgrupo}. 

In her thesis~\cite{melotese}, Melo first considered {\bf all} cyclic codes of $\F_q G$, that is, not only the minimals and 
computed dimension and minimum weights of these codes, using the following result.

\begin{lemma}{\em\cite[Proposição 2.1]{FDM}}\label{ttt}
Under the hypothesis above and considering $I_i$ the minimal ideal of $\F_q G$ generated by the primitive idempotent $e_i$, 
as in~\eqref{idempotPrimitivos}, for  $1\leq i\leq n$, we have
$$
d(I_i)=2|G_i|=2p^{n-i} \quad\mbox{ and }\quad \operatorname{dim}_{\F_q}I_i=\phi(p^i)=p^i-p^{i-1},
$$
and a basis for $I_i$ is 
$$
\mathcal{B}_i=\{ a(1-b)\widehat{G_i}\,|\, a \in A, 1\neq b \in B \},
$$
with $A$ a transversal of $G_{i}$ in $G_{i-1}$ and $B$ a transversal of $G_i$ in $G$.
For the minimal code $I_0=(\F_q G) e_0$, we have 
$$
w(I_0)=p^n \quad\mbox{ and }\quad \operatorname{dim}_{\F_q}I_0= 1.
$$
\end{lemma}

Considering that the dimension of a direct sum of ideals is the sum of their dimensions, Melo~\cite{melotese, meloart} focused her attention on computing minimum weight of the direct sum of minimal ideals as follows.

\begin{theo} \label{pesosFernanda}
Under the hypothesis of this section and of Lemma~\ref{ttt}, we have:

(i) {\em\cite[Lema 2.3]{melotese}} if $0< i < j$, then $w(I_i\oplus I_j) = 2|G_j|=2p^{n-j}$.

(ii) {\em \cite[Lema 2.4]{melotese} } If $1< j$, then $w(I_0\oplus I_j) = 2|G_j|=2p^{n-j}$.  

(iii) {\em \cite[Lema 2.5]{melotese}} If $I=I_0\oplus I_1$, then $w(I) = |G_1|=p^{n-1}$.

(iv)  {\em\cite[Lema 2.6]{melotese}} If $I=\bigoplus_{i=0}^{t} (\F_q G) e_i$, then $I = (\F_q G) \widehat{G_t}$ and 
$w(I) = |G_t|=p^{n-t}$.

(v)  {\em\cite[Lema 2.7]{melotese}} If $I=\bigoplus_{k=0}^{t} (\F_q G) e_{i_k}$, with 
$0\leq i_1 < i_2 < \cdots < i_t$ and  $ e_{i_1}+e_{i_2} + \cdots + e_{i_t}\neq  e_{0}+e_{1} + \cdots + e_{t} $, 
then $w(I) = 2|G_{i_t}|=2p^{n-i_t}$.
\end{theo}

Melo~\cite[Section 2.3]{melotese} also considered the distribution of weights for these cyclic codes. 
Furthermore, in~\cite[Chapter 3]{melotese}, she briefly compared cyclic and non-cyclic abelian codes 
of length $p^2$, fully exploring some examples using GAP Wedderga package. 

For the group $G=C_p\times C_p=\langle a\rangle\times\langle b\rangle$ and $\F_q$ a finite field
of $q$ elements such that $\bar{q}$ generates $U(\mathbb{Z}_p)$, the idempotents of $\F_q G$ are
$$
e_0=\widehat{G}, e_1=\widehat{a}-\widehat{G}, e_2=\widehat{b}-\widehat{G},f_i=\widehat{ab^i}-\widehat{G},
\mbox{ with } 1\leq i\leq p-1.
$$

Note that if $H$ and $K$ are any among the subgroups $\langle a\rangle$, $\langle b\rangle$,$\langle ab^i\rangle$, with  $1\leq i\leq p-1$,
then $G=H\times K$. For the idempotents $e=\widehat{H}-\widehat{G}$ and $e=\widehat{K}-\widehat{G}$ associated to $H$ and $K$, respectively, and considering the ideal $I= (\F_q G) e \oplus (\F_q G) f$, Melo proved:

\begin{theo} {\em\cite[Teorema 3.2.1]{melotese}}
The minimum weight of the ideal $I$ is $d(I)=2p-2$ and its dimension is $\operatorname{dim}_{\F_q} I = 2p-2$.
\end{theo}
\subsection{Essential idempotents an one weight cyclic codes}

In~\cite{CFM}, a special type of idempotent elements in
the semisimple group algebra of a finite abelian group is considered, the so called {\it essencial} idempotents.
These idempotents were previously considered by
Bakshi, Raka and Sharma in~\cite{barasha}, where they were called {\it non-degenerate}, in
the special case of group algebras of cyclic groups over finite fields.

\begin{definition}
In a semisimple group algebra $\F_q G$ of a finite group $G$, a primitive idempotent $e$ is an
{\bf essential idempotent} if $e\widehat{H}=0$, for every subgroup $H \neq \{1\}$ in $G$.
A minimal ideal of $\F_q G$ is called an {\bf essential ideal} if it is generated by
an essential idempotent.
\end{definition}
The following is a characterization of essential idempotents.
\begin{prop}{\em\cite[Proposition 2.3]{CFM}}
 Let $e\in \F_q G$ be a primitive central idempotent. Then $ e$ is
essential if and only if the map $\pi : G \longrightarrow Ge$ is a group isomorphism.
\end{prop}

\begin{cor}{\em\cite[Corollary 2.4]{CFM}}
If $G$ is an abelian  group and $\F_q G$ contains an essential idempotent,
then $G$ is cyclic.
\end{cor}

For cyclic groups, Chalom, Ferraz and Polcino Milies~\cite{CFM} proved the existence of a non-zero central 
idempotent which is the sum of all essential idempotents. They also give a criteria to determine essential
idempotents using the well-known Galois descent method  and, as a consequence, compute the
number of these idempotents in $\F_q C_n$, for $C_n$ a cyclic group of order $n$.

In~\cite[Section 3]{CFM} they  show that the coeficients of the primitive idempotents 
of a semisimple group algebra $\F_q A$, for $A$ is a finite abelian group,
can be easily computed as a concatenation of the coeficients of an essential
idempotent in the group algebras of a cyclic factor of $A$. In terms of coding
theory, this will imply that every minimal abelian code generated by a non
essential idempotent is a repetition code: their elements can be written as
repetitions of the coeficients of elements in a cyclic code generated by an
essential idempotent. In particular, one application of this is 
 to determine the
weight distribution of all codes when the weight distributions of codes
generated by essential idempotents are known.

Nascimento, in her Ph.D. Thesis~\cite{ruthtese}, uses this notion of essential idempotents 
to state conditions for a cyclic code in $\F_q C_n$ to be a one-weight code. 
Besides, she describes precisely the form of the elements on such a code and determines 
the number of one-weight codes in  $\F_q C_n$. She also constructs examples of two weight codes in  $\F_q (C_n\times C_n)$ and gives conditions to ensure that a code is of constant weight in  $\F_q A$, for $A$ an abelian group. 
Her work simplifies many of the proofs given by Vega~\cite{vega} for the same facts. In the literature there is also 
an interesting paper by Wood~\cite{woodconstant} on linear codes of constant weight.

\section{Codes over rings}\label{sec7}

In the 1990's many papers on cyclic codes over rings started to appear, motivated by the fact that good non linear binary codes were related to linear codes over $\mathbb{Z}_4$ (see, for example,~\cite{CS, KLP, PQ}). 
The paper~\cite{HKCSS} by Hammons et al. was even the best paper award for Information Theory of the IEEE-IT Society in the 1996 Symposium of IT - Whistler (Canadá). Wood~\cite{wood} addressed the problem of duality for modules over finite chain rings and applied it to equivalence of codes and to the extension theorem of MacWilliams.

In~\cite{CSlo} Carlderbank e Sloane determine the  structure of cyclic codes over $\mathbb{Z}_{p^m}$. 
Later on, in~\cite{KLP} Kanwar and López-Permouth did the same, but with different proofs. With the same techniques,
Wan~\cite{wan} extended the results from~\cite{KLP} to cyclic codes over Galois rings. 
Em 1999, Norton and S\u{a}l\u{a}gean-Mandache in~\cite{NSM} extended results of~\cite{CSlo, KLP} to 
cyclic codes over finite chain rings and later on, in 2004, 
Dinh and López-Permouth in~\cite{DLP} prove the same results in a different way. 

Codes over rings developed even more in the beginning of the 21st century that they deserved a 
CIMPA Summer School in 2008~\cite{sole}. Further works can be found in~\cite{Dinh},~\cite{Liu},~\cite{MMR}. 
A small survey on the subject is~\cite{greferath}.

In his thesis~\cite{AndTese,AndArt}, Silva used group ring approach to characterize cyclic codes over chain rings, their duals and some conditions on self-dual codes, simplifying the proofs and improving results given in~\cite{DLP}.

Let $R$ be a finite commutative chain ring with unity such that  $\mid R \mid=q^k$, for a prime $q$. 
For $M$ the maximal ideal of $R$, the quotient $\overline{R}=\frac{R}{M}$ is a field and
we work under the hypotesis that  $q\nmid \mid G \mid$, for a finite cyclic group $G$. 
Under these conditions, the group ring $RG$ is a principal ideal ring, as Silva proves in~\cite[Teorema 2.1.9]{AndTese},
after characterizing all the ideals in $RG$. The following general fact is a basis for all this work.

\begin{theo}{\em\cite[Teorema 2.1.2]{AndTese}}
Let $R$ be a local ring, with  maximal ideal  $M=<a>$ and $\mid R \mid=q^k$,  and $G$ a cyclic group of order  $n$ such that $q \nmid n$. If
$\{\overline{e}_0,...,\overline{e}_m\}$ is a full set of primitive orthogonal idempotents in  $\overline{R}G$, then 
$\{e_0,...,e_m\}$ is a full set of primitive orthogonal idempotents in $RG$.
\end{theo}

The next theorem characterizes all cyclic codes of length $n$ over the local ring $RG e_i$ (see~\cite[Corollary 11.31]{AndTese}), for $R$ a chain ring and $e_i$ a primitive orthogonal idempotent, translating results of~\cite{DLP} to the
group ring setting. 
To simplify the notation we write  $(RG)a^je_i$ as $\left<a^je_i\right>$.

\begin{theo}{\em \cite[Teorema 2.1.3]{AndTese}}{\label{car}}
Let $R$ be a commutative finite chain ring with unity, $\mid R \mid=q^k$, $M=\left<a\right>$ the maximal ideal of $R$ and $t$  the nilpotency index o índice de nilpotência of $a$ in $R$.  Let  $G=C_n$ such that $q\nmid n$. If $I$  is an ideal of 
$RGe_i$, then $I$ is of the form $I=\left<a^{k_i} e_i\right>$, with $0\leq k_i \leq t$. 
\end{theo}

\begin{cor}{\em \cite[Corollary 2.1.4]{AndTese}} Under the same hypothesis of Theorem~\ref{car}, the ideal  $RGe_i$ is 
indecomposable in $RG$ and the code  $\left<a^{t-1}e_i\right>$ is minimal.
\end{cor}

From this we have a characterization of all cyclic codes of length $n$ over chain rings.

\begin{theo}\label{caracterizacao}
Let $R$ be a commutative finite chain ring with unity, $\mid R \mid=q^k$, $M=\left<a\right>$ the maximal ideal of $R$ and $t$  the nilpotency index of $a$ in $R$.  Let  $G=\left<g_0\,/\,g_0^n=1\right>$ be such that $q\nmid n$ and
$\{e_0,...,e_m\}$ be a full set of primitive orthogonal idempotents of $RG$. Then:

(i) {\em\cite[Teorema 2.1.5]{AndTese}} If $I$  is an ideal of 
$RG$, then $I$ is of the form
 $I=I_0\oplus...\oplus I_m$, with $I_i=\left<a^{k_i}
e_i\right>$, for $0\leq k_i \leq t$.

(ii) {\em \cite[Teorema 2.1.8]{AndTese}} The number of such codes of length $n$ over  $R$ is $(t+1)^{m+1}$.
\end{theo}

One important data in a code is its number of words. Next theorem gives this number for cyclic codes over finite chain rings.
We have 

  $$RG=RGe_0\oplus...\oplus RGe_m\simeq \frac{R[x]}{\langle x^n-1\rangle}\simeq \frac{R[x]}{\langle f_0 \rangle}\oplus...\oplus \frac{R[x]}{\langle f_m \rangle},$$
where $f_i$ are irreducible factors of $x^n-1$ and , after reordering the indexes if necessary, we have 
$RGe_i \simeq \frac{R[x]}{\langle f_i \rangle}$. Hence, $\mid RGe_i \mid=\mid R \mid^{w_i}$,
for $w_i=deg(f_i)$. 

\begin{theo}{\em\cite[Teorema 2.1.7]{AndTese}}\label{conta} Under the same hypothesis of Theorem~\ref{car}, 
let $C$ be a cyclic code of the form $C=\left<a^{k_{i_1}}e_{i_1}\right>\oplus...\oplus \left<a^{k_{i_r}}e_{i_r}\right>$ in $RG$. The the number of words in $C$
is $\mid C \mid= \mid \overline{R} \mid^{\displaystyle{\sum_{s=1}^{r}}(t-k_{i_s})w_{i_s}}$.
\end{theo}

Considering $*:RG\longrightarrow RG$ the classical involution, Silva also gives a description of the dual cyclic codes 
in $RG$ as follows.

\begin{theo}{\em\cite[Teorema 2.2.3]{AndTese}}\label{dual code}
 Under the same hypothesis of Theorem~\ref{caracterizacao}, the dual code of
a cyclic code $C=\left< a^{k_0}e_{0}\right>\oplus...\oplus\left<a^{k_{m}}e_{m}\right>$, with $0\leq k_i\leq t$, is 
$C^{\perp}=\oplus \sum_{r=0}^{m}\left<a^{t-k_{r}}{e_{r}}^*\right>$.
\end{theo}

As in~\cite{DLP}, Silva in~\cite[Section 2.2]{AndTese} states the conditions for the ring $R$ under which the 
group ring $RG$ admits self-dual codes.

Chapter 3 of~\cite{AndTese} is dedicated to codes over chain rings of length $p^n$, for a prime $p$, 
extending the results of Ferraz and Milies~\cite{FM} and of Melo~\cite{melotese} to this context. 
Silva also proves in~\cite[Teorema 3.0.14]{AndTese} some facts 
about the size of such codes and computes minimum weight of these codes
~\cite[Teoremas 3.0.15 to 3.0.18]{AndTese}, similarly to Theorem~\ref{pesosFernanda}. He also discusses about free codes 
in $RG$ in~\cite[Section 3.1]{AndTese} and about MDS codes of length $p^n$ over $R$ in~\cite[Section 3.2]{AndTese}.
Finally, in~\cite[Chapter 4]{AndTese}, Silva proves all such results for cyclic codes of length $2p^n$ over finite chain rings.

There are also interesting discussion on equivalence of linear codes over rings in~\cite{DLP2, DLP3, wardwood,  wood2}.

% fim do small para as letras
\end{document}